\documentclass[11pt,a4paper,openright]{article}
\usepackage[affil-it]{authblk}
\usepackage{amsmath}
\usepackage[numbers]{natbib}
\usepackage{amsfonts}
\usepackage{amssymb}
\usepackage{mathtools}
\usepackage{mathrsfs,bbm}
\usepackage{xcolor}
\usepackage{lmodern}
\usepackage{hyperref}

\usepackage{graphicx}
\usepackage{epstopdf}

\usepackage{bm}
\usepackage{accents}
\usepackage{orcidlink}
\usepackage{cleveref}
\usepackage[ruled,vlined]{algorithm2e}

\newcommand{\Rbb}{\mathbb{R}} 
\DeclareMathOperator{\Mat}{Mat}

\DeclareMathOperator{\orth}{orth}

\DeclareMathOperator{\St}{St}
\newcommand{\bk}{\mathbf{k}}
\newcommand{\x}{\mathbf{x}}

\newcommand{\dist}{{\operatorname{\tt spec\_dist}}}

\begin{document}

\title{A Low-Rank BUG Method for Sylvester-Type Equations}
\author[1]{Georgios Vretinaris\orcidlink{0009-0007-6819-1586}}

\affil[1]{Mathematisches Institut, Universit{\"a}t T{\"u}bingen, Auf der Morgenstelle 10, 72076 T{\"u}bingen, Germany.\thanks{ \texttt{ georgios.vretinaris@uni-tuebingen.de}}}
\date{}

\maketitle

\begin{abstract}
	We introduce a low-rank algorithm inspired by the Basis-Update and Galerkin (BUG) integrator to efficiently approximate solutions to Sylvester-type equations. The algorithm can exploit both the low-rank structure of the solution as well as any sparsity present to reduce computational complexity. In the worst case scenario, i.e. even when a standard dense solver, such as the Bartels-Stewart algorithm, is used for the reduced Sylvester equations generated by our approach, the overall computational complexity for constructing and solving the associated linear systems reduces to $\mathcal{O}\left(kr(n^2+m^2 +mn + r^2)\right)$, for $X\in \Rbb^{m \times n}$, where $k$ is the number of iterations and $r$ the rank of the approximation, but any method can be used to exploit existing structure to optimize the solution even further.
\end{abstract}

\section{Introduction}

The Sylvester equation
\begin{equation}
	A X + X B^T = C, \label{eq:sylv}
\end{equation}
where $A\in \Rbb^{m \times m},\, B \in \Rbb^{n \times n}$ and $C \in \Rbb^{m \times n}$, is a classical linear matrix equation that arises in various different areas of applied mathematics, most commonly in control theory \cite{Datta2003-ej}, model reduction \cite{Antoulas2005-yj} and numerics of partial differential equations.
In these settings $A$ and $B$ are often large, sparse, or structured matrices. Efficiently solving this equation is therefore a fundamental task. Of course, the solution of Eq. \eqref{eq:sylv} exists uniquely if and only if the spectra of $A$ and $-B$ are disjoint. For further reading on the equation, the reader is advised to look into \cite{Bhatia1996-tg,Bhatia1997-wn,Higham2002-av}.

By Penzl \cite{Penzl2000}, we are already aware of the fact that in the case of the continuous-time algebraic Lyapunov equation (when $B = A$) with symmetric and c-stable coefficient matrix $A \in \Rbb^{n\times n}$, i.e. it's spectrum lies in the left-half of the complex plane, and a low-rank symmetric right-hand side matrix $C \in \Rbb^{n\times n}$, the solution $X \in \Rbb^{n\times n}$ shows a fast decay of the nonincreasingly ordered eigenvalues. There are other works which extend this endeavor \cite{Baker2015,Antoulas2002,sabino2006lyapunov}.
Beckermann and Townsend proved that in the case of normal coefficient matrices, the singular values show a decay which depends on the distance between the spectra of $A$ and $-B$ \cite{Beckermann2017}. This result was recently augmented by moderately relaxing the normality condition \cite{Clouatre2025-kd}.
The consequence of these findings is that we can now know when a low-rank approximation is justified.

The importance of the Sylvester (and Lyapunov) equation in many applications has created an extensive amount of literature dedicated to efficiently solving for a low-rank approximation of its solution. Most algorithms tend to fall be in some large class of methods, for example, the low-rank alternate direction implicit (LR-ADI) methods \cite{Li2002-hl, Benner2009-wo}, Krylov subspace methods \cite{Casulli2024-vd,Casulli2025,Palitta2025-og}, mixed-precision methods \cite{Benner2025-ww,Schulze2025-iv} and methods that are part of more than just one class. The reader can find a more extensive survey in \cite{Simoncini2016-ds}

We present a novel low-rank approach to solving the Sylvester equation, inspired by the Basis-Update and Galerkin (BUG) integrator recently introduced for dynamical low-rank approximation \cite{Ceruti2022-ty}. The idea of the BUG method is essentially to split a large matrix problem into three smaller ones. This way, both the storage and the number of operations can be significantly reduced. Unfortunately, a rigorous convergence proof of the proposed algorithm is still missing.

This work is structured as follows: In Section 2 we will present the algorithm for solving the Sylvester Eq. \eqref{eq:sylv} and take the reader through the main ideas of the BUG method. Section 3 generalizes the algorithm to solving the multilinear or tensor Sylvester equation and provides some comments on how to avoid creation of dense matrices. Numerical experiments will be presented in Section 4.

\section{Low-Rank BUG Sylvester Solver}

If a rank-$r$ solution $X\in \Rbb^{m\times n}$ of Eq. \eqref{eq:sylv} exists, it can always be expressed in the form $X = USV^T$, where $U \in \Rbb^{m \times r}, V\in \Rbb^{n\times r}$ are thin matrices with orthonormal columns, or in other words elements of Stiefel manifolds $\St(n,r) \coloneqq \{A \in \Rbb^{n\times r}\, \colon\, A^T A = I_r\}$ and $S \in \Rbb^{r \times r}$ is nonsingular. For the singular value decomposition (SVD), $S$ is a diagonal matrix containing the singular values in descending order, but we do not need such special form here. Using this fact in Eq. \eqref{eq:sylv} would yield the following expression
\begin{equation}
	 A USV^T + USV^T B^T = C \label{eq:usv}
\end{equation}
The idea is to ``hide" the $S$ matrix and focus on the bases. We therefore choose to abbreviate $US\eqqcolon K$ and create a new Sylvester equation for $K$ of lower dimension $n\times r$, by post-multiplying both sides of Eq. \eqref{eq:usv} by $V$, which now writes
\begin{equation}
	 A K + K(V^T B V)^T = C V. \tag{K-step}
\end{equation}
Similarly, for the right basis $V$, define $L \coloneqq VS^T$ and multiply by $U$ to get
\begin{equation}
	B L + L (U^T A U)^T = C^T U. \tag{L-step}
\end{equation}
Although $K$ and $L$ do not directly recover $X$, they suffice to extract the bases via orthogonalization, for example by a QR decomposition, (which we denote as the operation $\orth$), and get
\begin{equation}
\hat{U} = \orth(K), \qquad \hat{V}=\orth(L).
\end{equation}
Given the new bases, we can accurately solve for the $S$ matrix provided by the Galerkin condition, which now is a heavily reduced Sylvester equation of dimension $r \times r$:
\begin{equation}
	 (\hat{U}^T A \hat{U}) S + S (\hat{V}^T B \hat{V})^T = \hat{U}^T C\hat{V} \tag{S-step}.
\end{equation}
If we knew the exact $U$ and $V$ a priori, solving for $K$ and $L$ would just yield $\hat{U}=U$ and $\hat{V} = V$. Since the bases are unknown, we solve for them iteratively starting from initial points which can be chosen randomly since they lie in Stiefel manifolds, which are compact.

To recap, the fundamental idea is that one should update the left and right low-rank bases the equation onto reduced subspaces and once the correct bases are found one can move on to compute the ``singular value" matrix by the Galerkin condition. In practice this translates to splitting the initial Sylvester equation into three reduced ones.
The discussion above immediately motivates our BUG inspired method for the Sylvester equation. Namely:

\begin{algorithm}[H]
\caption{Fixed-rank BUG Sylvester Solver}
\label{alg:fixed-rank}
\SetAlgoLined
\DontPrintSemicolon
\KwIn{Coefficient matrices $A \in \Rbb^{m\times m},B \in \Rbb^{n \times n}$, initial guesses $U_0\in \St(m,r),\, V_0\in \St(n,r)$, abs. tolerance {\tt tol}}
\KwResult{Matrix $X = U S V^T$}
\Begin{
Initialize {\tt prev\_err} = $10^{10}$ \\
\For{$i = 0$ to {\tt max\_iter}}
{
    Solve for $K_{i+1}$, $A K_{i+1} + K_{i+1} V_i^T B^T V_i = C V_i$  \tcp*{K-Step}
    Get $U_{i+1} = \orth(K_{i+1})$\\
    Solve for $L_{i+1}$, $B L_{i+1} + L_{i+1} U_i^T A^T U_i = C^T U_i$     \tcp*{L-Step}
    Get $V_{i+1} = \orth(L_{i+1})$\\
    Solve for $S_{i+1}$, $U_{i+1}^T A U_{i+1} S_{i+1} + S_{i+1} V_{i+1}^T B V_{i+1} = U_{i+1}^T C V_{i+1}$    \tcp*{S-Step}
    Define $X_{i+1} \coloneqq U_{i+1} S_{i+1} V_{i+1}^T$ \\
    Define {\tt curr\_err} = $\|A X_{i+1} + X_{i+1}B^T - C\|_F $\\
    \If{{\tt curr\_err} $\leq {\tt tol}$ {\bf or} $|{\tt curr\_err} - {\tt prev\_err}| \leq {\tt tol}$}{
        \textbf{break}
    }
	Update {\tt prev\_err} = {\tt curr\_err} 
}
}
\end{algorithm}

Notice that the fixed-rank algorithm unfortunately requires taking a guess of the correct rank, which might not be the most optimal one. This is easily remedied by following the tactics of \cite{Ceruti2022-ig} which provide a rank-adaptive method given a tolerance.

\begin{algorithm}[H]
\SetAlgoLined
\DontPrintSemicolon
\caption{Rank-adaptive BUG Sylvester Solver}
\label{alg:cap}
\KwIn{Coefficient matrices $A \in \Rbb^{m\times m},B \in \Rbb^{n \times n}$, initial guesses $U_0\in \St(m,r),\, V_0\in \St(n,r)$, abs. tolerance {\tt tol}, $c$ stagnation weight, truncation tolerance $\vartheta$}
\KwResult{Matrix $X = USV^T$}
\Begin{
Initialize {\tt prev\_err} = $10^{10}$ \\
\For{$i = 0$ to {\tt max\_iter}}{
    Solve $A K_{i+1} + K_{i+1} V_i^T B^T V_i = C V_i$ \tcp*{K-Step}
    Get $\widehat{U}_{i+1} = \orth([K_{i+1}\quad U_i])$    \\
    Solve $B L_{i+1} + L_{i+1} U_i^T A^T U_i = C^T U_i$ \tcp*{L-Step}
    Get $\widehat{V}_{i+1} = \orth([L_{i+1}\quad V_i])$  \\
    Solve $\widehat{U}_{i+1}^T A \widehat{U}_{i+1} S + S \widehat{V}_{i+1}^T B \widehat{V}_{i+1} = \widehat{U}_{i+1}^T C \widehat{V}_{i+1}$ \tcp*{S-Step}
    Decompose $S = P \Sigma Q^T$\\
    Set $\hat{r} \coloneqq \min\limits_{r}\left| \left(\sum\limits_{j=r+1}^{2r_i} \sigma_{j}^2 \right)^{1/2} - \vartheta \right|$\\
    Truncate up to $\hat{r}$\\
    Define $U_{i+1} = \widehat{U}_{i+1}P_{\hat{r}}$, $V_{i+1} = \widehat{V}_{i+1}Q_{\hat{r}}$, $S_{i+1}=\Sigma_{\hat{r}}$\\
    Define $X_{i+1} \coloneqq U_{i+1} S_{i+1} V_{i+1}^T$ \\
    Define {\tt curr\_err} = $\|A X_{i+1} + X_{i+1}B^T - C\|_F $\\
	\If{{\tt curr\_err} $\leq {\tt tol}$ {\bf or} $|{\tt curr\_err} - {\tt prev\_err}| \leq c\cdot{\tt tol}$}{
		\textbf{break}
	}
Update {\tt prev\_err} = {\tt curr\_err} 
}
}
\end{algorithm}

The crux of the iteration is the update of the bases. In fact, the S-step can be removed from the loop and performed once after convergence of the bases is achieved. Once again, a formal convergence analysis is still missing, so a first naive idea is to perform the S-step each time, as it is also the most inexpensive system to solve and check whether the residual is below the desired tolerance. Furthermore, since neither the solution nor the rank is known a priori, there is tension between the stopping criterion and the truncation tolerance which is why a stagnation clause is added to the stopping criterion.

The standard algorithm for the direct solution of eq. \eqref{eq:sylv} for full matrices is the Bartels-Stewart algorithm \cite{Bartels1972-jh}, which requires $\mathcal{O}(n^3+m^3)$ operations for the Schur decompositions and $\mathcal{O}(nm^2+n^2m)$ to get the reduced linear systems. In the worst case scenario, i.e. when the coefficient matrices are dense and unstructured, if the K-, L- and S-steps are solved using this algorithm then the complexity of the Schur decompositions is still the same (though they only need to happen once for each coefficient matrix), but the operations for the creation and solution of the linear system drops to $\mathcal{O}\left(kr(n^2+m^2 +mn + r^2)\right)$, where $k$ is the number of iterations needed and $r$ the final rank. On the other hand, if the coefficient matrices $A$ and $B$ are sparse, or/and symmetric, as often occurs in applications, one can use methods that exploit their sparsity and substantially reduce the number of operations needed. In any case, one can use the most optimal method available, in order to solve each of the steps.

\section{Extension to Tucker Decompositions} \label{sec:extension}

Since low-rank methods also aim to help remedy the curse of dimensionality, a natural next step is to generalize the BUG approach to multidimensional tensors represented in Tucker format. Before doing so, we have to introduce some notation. Let $X \in \Rbb^{I_1,\times I_2\times \ldots \times I_N}$ and $A \in \Rbb^{J\times I_n}$, then the n-mode product is written $X \times_n A$, is of size $I_1,\times \ldots \times I_{n-1} \times J \times I_{n+1}\times \ldots \times I_N$ and elementwise is represented in the following way:

\[(X \times_n U)_{i_1\ldots i_{n-1}j i_{n+1} \ldots i_N} = \sum_{{i_n}=1}^{I_n} x_{i_1 i_2 \ldots i_N}u_{ji_n}.\]

Finally, we define the matricization of a tensor $X$ as before, as $\Mat_n{X} \in \Rbb^{I_n \times I_{\neg n}}$, where $I_{\neg n} = \prod_{j\neq n} I_j$. We recall that the $n$-th matricization aligns in the $k$-th row (for $k=1,\ldots, I_n)$ all entries of $X$ that have the index $k$ in the $i$-th position, usually ordered co-lexicographically. For further reading and examples, the reader is advised to continue with \cite{Kolda2009-km,Ceruti2021-od}.

Let us consider the following common scenario in numerics of differential equations, when discretizing a Laplacian in higher dimensions, yielding a tensor Sylvester equation:

\begin{equation}
    \sum_{i=1}^d X \times_i A_i = B. \label{eq:tensor-sylv}
\end{equation}

Let $X\in \Rbb^{n_1 \times \ldots \times n_d}$ of low-(multilinear-)rank, we can then write it via its Tucker Decomposition
\begin{equation}
	X = C \times_1 U_1 \times_2 U_2 \times_3 \ldots \times_d U_d = C \bigtimes_{i=1}^{d} U_i,
\end{equation}
where the size of $C$ is $r_1 \times \ldots r_d$, and $U_i \in \Rbb{n_i \times r_i}$. Using this form of $X$ in Eq. \eqref{eq:tensor-sylv}, casts it into
\begin{align*}
    &\sum_{i=1}^d C \times_1 U_1 \times_2\ldots \times_i (A_i U_i) \times_{i+1} \ldots \times_{d-1}U_{d-1} \times_d U_d = B   .
\end{align*}

Performing a matricization of this Tucker tensor over mode $k$, in the general case, reads

\begin{equation}
	\begin{aligned}
		    A_k U_k \Mat_k(C) \bigotimes_{j\neq k} U_j^T+ \sum_{i\neq k}U_k \Mat_k(C) \left( U_d^T \otimes \ldots \otimes(A_i U_i)^T \otimes \ldots \otimes U_1^T \right)\\
		    = \Mat_k (B) 
	\end{aligned}
\label{eq:Tucker1}
\end{equation}

where we used the formula:
\begin{equation}
	\Mat_i\left(C \bigtimes_{i=1}^d U_i\right) = U_i \Mat_i(C) \left(\bigotimes_{j\neq i}U_i^T \right)
\end{equation}

We notice that if we write $\left(\Mat_k(C)\right)^T = Q_k S_k^T$ and then define the $K_k \coloneqq U_k S_k$, $V_k^T \coloneqq Q_k^T \bigotimes_{j\neq k} U_j^T$
and 

\[\widetilde{U}_{ij} \coloneqq \begin{cases}
    U_j  &\text{if} \,\, j\neq i, \\
    A_i U_i &\text{if} \,\, j = i
\end{cases}\]

Eq. \eqref{eq:Tucker1} is rewritten into

\begin{align}
    A_k K_k + K_k Q_k^T \sum_{i\neq k}^d   \bigotimes_{j\neq k} \widetilde{U}_{ij}^T V_k = \Mat_k(B) V_k \\
    A_k K_k + K_k \sum_{i\neq k}^d P_{ik} = \Mat_k(B) V_k,
\end{align}
where we also multiplied both sides by $V_k$.

Thus, we ended up with a Sylvester equation for the basis matrix $U_k$ which we know how to solve. The matrix $P_{ik}$ can be compactly expressed as follows

\begin{align*}
    P_{ik}&\coloneqq\left(Q_k^T \bigotimes_{j\neq k} \widetilde{U}_{ij}^T\right)V_k  \\
    &= \left(Q_k^T \bigotimes_{j\neq k} \widetilde{U}_{ij}^T\right) \left( \left(\bigotimes_{j\neq k} {U}_j\right)  Q_k\right) \\
    &= Q_k^T \cdot [I_{r_d} \otimes \ldots \otimes U_i^T A_i^T U_i \otimes \ldots \otimes I_{r_1}] \cdot Q_k \\
    &= Q_k^T R_i Q_k.
\end{align*}

%

To deal with the core tensor, we multiply both sides of Eq. \eqref{eq:tensor-sylv} by $\bigotimes_{i=1}^d U_i^T$ which writes:
\begin{align}
	&\sum_{i=1}^d C \times_i M_i^T = B \bigtimes_{i=1}^d U_i^T, \label{eq:tucker-core}\\
     \text{or,} \qquad&\Mat_0(C) \sum_{i=1}^d R_i = \Mat_0(B) \bigotimes_{i=1}^d U_i. \label{eq:tucker-lin}
\end{align}
where $M_i \coloneqq U_i^T A_i^T U_i$. Note that the $R_i$'s are matrices of size $\tilde{r} \times \tilde{r}$, with $\prod_{i=1}^d r_i$, which unfortunately still scales exponentially with the dimension of the problem. Exploring further hierarchical decompositions like the Tree Tensor Networks (TTN) can help solve this problem. But the multilinear system is now significantly smaller than the initial one, allowing direct solutions to the tensor Sylvester equation for the core.

With the aforementioned equations in mind, the algorithm that solves the tensor Sylvester equation using Tucker Tensors emerges quite naturally in the following way.

\begin{algorithm}[H]
\SetAlgoLined
\DontPrintSemicolon
\caption{Fixed-rank BUG Tensor Sylvester Solver}
\KwIn{Coefficient matrices $(A_k \in \Rbb^{n_k\times n_k})_{k=1}^d$, initial guesses $\left(U_k^0\in \St(n_k,r)\right)_{k=1}^d$, $C\in \Rbb^{n\times\ldots\times n}$, abs. tolerance {\tt tol}}
\For{$i = 0$ to {\tt max\_iter}}{
Initialize {\tt prev\_err} = $10^{10}$ \\
    \For{$k = 1$ to $d$, {\rm in parallel},}{
        Compute the QR decomposition $\phantom{x}\qquad\Mat_k(C)^T = Q_k S_k^T \in \Rbb^{r_{\neg i} \times r_i}$

        Solve $A_k K_k^{i+1} + K_k^{i+1} \sum_{j\neq k}^d P_{jk}^i = \Mat_k(B) V_k^i$
        
        Define $U^{i+1}_k \coloneqq \orth(K_k^{i+1})$
    }
    Solve $\sum_{i=1}^d C \times_i M_i^T = B \bigtimes_{i=1}^d U_i^T$
    
    Define $X^{i+1} \coloneqq {\rm Ten}_0(C^{i+1}) \bigtimes_{j=1}^d U_j^{i+1}$

    Define {\tt curr\_err} =  $\|\sum_{j=1}^d X^{i+1} \times_j A_j - B\|_F$\\
    	\If{{\tt curr\_err} $\leq {\tt tol}$ {\bf or} $|{\tt curr\_err} - {\tt prev\_err}| \leq {\tt tol}$}{
    	\textbf{break}
    }
    Update {\tt prev\_err} = {\tt curr\_err} 
}
\end{algorithm}

We should underline here that he products $\Mat_k(B) V_k$ and $\Mat_0(B) \bigotimes_i U_i$ are never computed as dense matrices, as this would contradict the entire effort of using the low-rank structure. Since we assume $B$ to have a low-rank structure i.e. it should also be a Tucker Tensor in this scenario, one can end up in the linear system of Eq. \eqref{eq:tucker-lin} by first performing the matrix multiplications in each basis and then get a reduced basis which can n-mode multiplied to the core of $B$.

Furthermore, one can get the rank-adaptive method by augmenting the bases exactly as shown in the matrix case, and then doing a HOSVD in the Tucker case.

\section{Numerical Experiments}

In all of the following, the (multilinear-)rank of the right-hand side is chosen to be 7, and all coefficient matrices lie in $\Rbb^{n \times n}$, for $n=16384$, unless stated otherwise. Furthermore, we define $\mathcal{L}(X) \coloneqq \sum_{i=1}^d X \times_i A_i$ as the Sylvester operator, $Y$ to be the numerical low-rank approximation, $X_r$ to be the densely computed solution, truncated up to rank-$r$ and $R\coloneqq \mathcal{L}(Y) - B$ the residual, where $B$ is the right-hand side. Recall that we are using random initial guesses for the bases matrices $U_i$ and the same holds for the core tensors, since we only need the $Q$ factors from their QR decomposition. The ranks of the initial guesses are exactly the ranks of the right-hand side. \footnote{The source code that reproduces the figures presented can be found in \href{https://github.com/gvretina/LR_BUG_Sylvester.jl}{https://github.com/gvretina/LR\_BUG\_Sylvester.jl}}

For the tensor case, since the connection tensors lack the nice structure of the $\Sigma$ matrix in the SVD, the generalized ``singular-values" won't be shown. We opt for the three-dimensional scenario and we use $n=4096$. Otherwise everything remains the same.

\subsection{Poisson Equation (Sparse Structured Case)}

Since the primary motivation to develop this method came from solving a Poisson problem, we shall start by introducing it.

\begin{equation}
	\Delta u = f \Rightarrow \sum_{i=1}^d X \times_i D_i = B, 
\end{equation}

\noindent where $D_i \coloneqq ({\tt tridiag}(1,-2,1))/\Delta x_i^2 $, and 
\[\widetilde{B}\coloneqq \sum\limits_{\bk\in \{-3,\ldots,3\}^d} a_\bk \cos(\bk \cdot \x + \varphi_\bk),\] 

with $\x \in [0,4\pi]^d$, $a_\bk \coloneqq {\tt randn}()/(1+\|\bk\|_1)$, $\varphi_\bk \coloneqq 2\pi *{\tt rand}()$ and $d$ being the dimension ($d=2$ for the matrix case and $d=3$ for the Tucker Tensor case). 
\subsubsection{Matrix Case}

\begin{figure}[ht!]
	\centering
	\includegraphics[width=\textwidth]{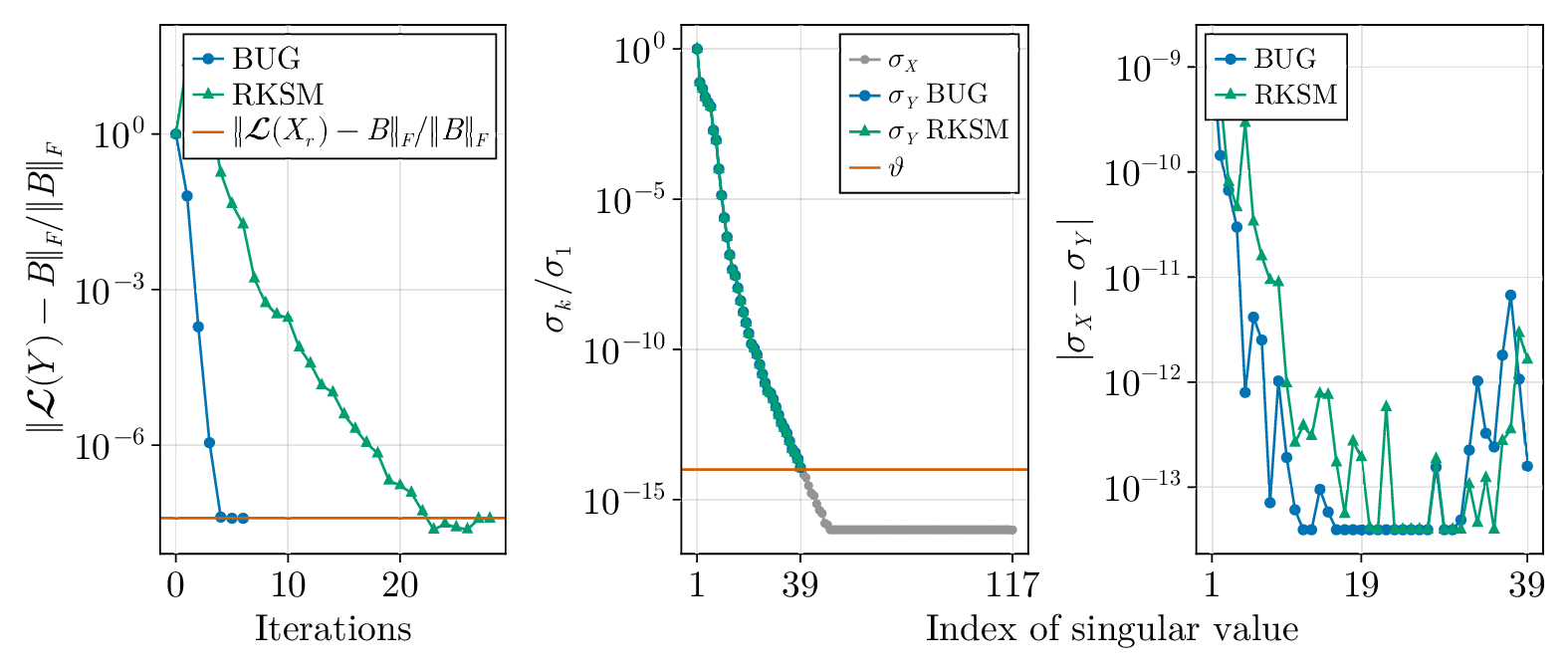}
	\caption{Numerical behavior for the Poisson equation with Dirichlet boundary conditions. {\bf Left figure}: Convergence of the rank-adaptive algorithm and the block rational Krylov method by \cite{Casulli2024-vd}. {\bf Middle figure}: The singular values of the solution $X$ given by a direct numerical solution, over its largest singular value, alongside the singular values of the low-rank approximation and the truncation tolerance $\vartheta = 10^{-14}\|\Sigma_Y\|_F$. {\bf Right figure}: The distance between the singular values of the exact solution to the low-rank approximation.}
	\label{fig:laplacian-dirichlet}
\end{figure}

In order to portray the power of this method, let's see how we deal with each case numerically too. Since the K- and L- steps are equivalent we will only show it for one of them. Without loss of generality, $AK + KM^T = CV$, now since $M = V^T D V$, it is obvious that $M$ is symmetric allowing for an $LDL^T$ decomposition and get $A\tilde{K} +\tilde{K} \Lambda = \tilde{Z}$ where $\tilde{K} \coloneqq KE$ and $\tilde{Z} =CVE$. Since $\Lambda$ is diagonal, $\tilde{K}$ and $\tilde{Z}$ are thin matrices we can solve this for each column and get $r$ tridiagonal systems. The S-step is small can can always just be solved with the standard Bartels-Stewart. As shown in \Cref{fig:laplacian-dirichlet}, the BUG method requires less than 10 iterations of the aforementioned scheme and stagnates due to the truncation tolerance policy, but identifies correctly the $r$ leading singular values.

\subsubsection{Tucker Tensor Case}

\begin{figure}[ht!]
	\centering
	\includegraphics[width=0.45\textwidth]{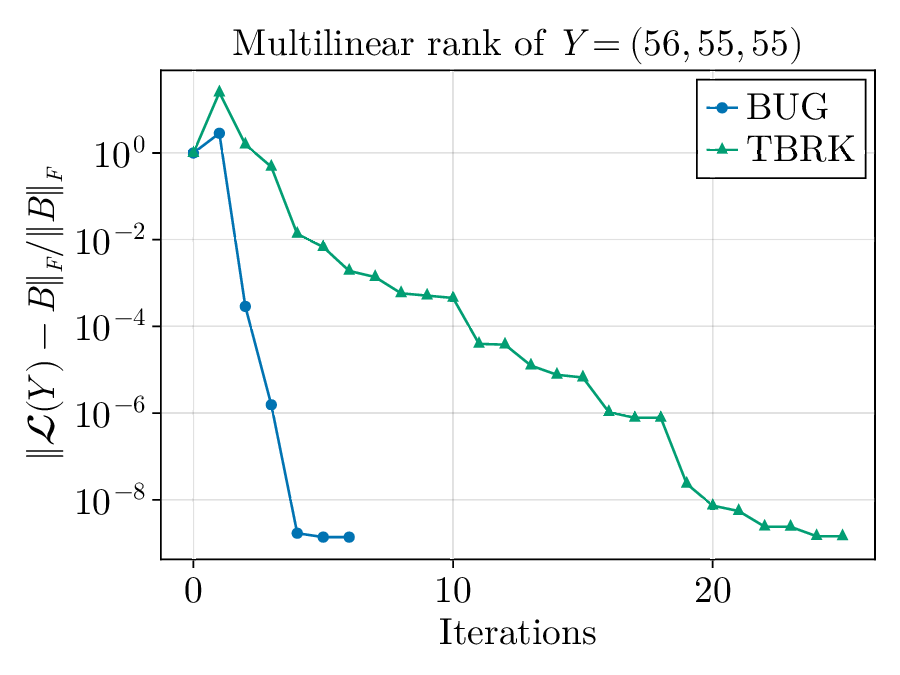}
	\caption{Convergence of the rank-adaptive method and the block rational Krylov method by \cite{Casulli2025} for Poisson problem with Dirichlet boundary conditions. In both cases, the truncation tolerance was set to $\vartheta = 10^{-14}\|C_Y\|_F$.}
	\label{fig:tucker-poisson}
\end{figure}

Everything we discussed about in the matrix case are of course directly applicable in the tensor case too. Which is further supported by \Cref{fig:tucker-poisson}. Of course, the computing the dense exact solution is not feasible now thus we are not able to present $X_r$ and the relative residual floor due to truncation.

It is obvious that the algorithm performs in much the same way when augmenting the dimensions and using Tucker Tensors instead of matrices.

\subsection{Random Matrices (Dense Unstructured Case)}

A good way to test the proposed method is with completely random matrices for which we only control the spectral distance of the coefficient matrices. This way we confirm our low-rank approximability conjecture from the introduction and we can test our methods in dense and unstructured scenaria. To do so, we construct $A_i \coloneqq P_i \Lambda_i {P_i}^{+}$, where $P_i = {\tt rand}(n,n)$, and normalize each column of $P_i$, ${P_i}^+$ is the typical Moore-Penrose inverse, $\Lambda_i = (-1)^i{\rm diag}({\tt rand}(n)+1+\dist(i-1))$. Finally $B \coloneqq C \bigtimes_{i=1}^d U_i$, where $C = {\tt rand}\underbrace{(7,\ldots,7)}_{d \text{ times}}$ and $U_i \in \St(n,r)$ arbitrarily chosen.

\subsubsection{Matrix Case}

\begin{figure}[ht!]
	\centering
	\includegraphics[width=0.95\linewidth]{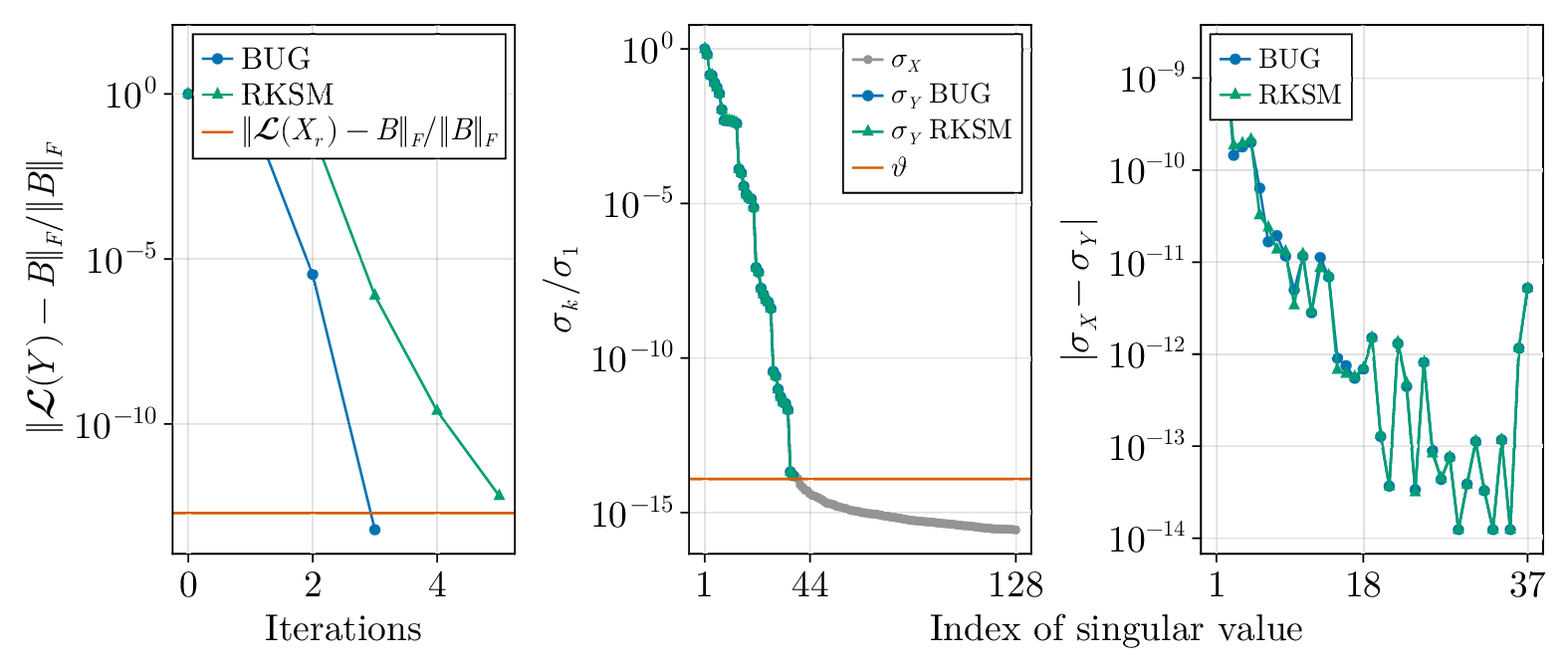}
	\caption{Numerical behavior for the random case. {\bf Left figure}: Convergence of the rank-adaptive algorithm and the block rational Krylov method by \cite{Casulli2024-vd} for matrices. {\bf Middle figure}: The singular values of the solution $X$ given by a direct numerical solution, over its largest singular value, alongside the singular values of the low-rank approximation and the truncation tolerance $\vartheta = 10^{-14}\|\Sigma_X\|_F$. {\bf Right figure}: The distance between the singular values of the exact solution to the low-rank approximation.}
	\label{fig:random}
\end{figure}

Since we have already seen the scenario Lyapunov equation in the previous experiment, which for the case of the Laplace operators also treats the small eigenvalue situation, we try to portray how increasing the spectral radius is enough to get a solution that exhibits low-rank structure. For the matrix case specifically, we choose the spectral distance $\dist=10$. Once again, we can exploit the structure of each substep of our method to end up with a faster algorithm overall. In this case, since we lack structure from $A$ or $M$ (for the K-step), we have to opt for a Krylov method since densely computing everything isn't feasible for large-scale simulations. But we don't need to build the Krylov space anew for each iteration. Let $\mathcal{K}_m(A,C_u)$ be the typical Krylov space, where $C = C_u C_s C_v^T$. Building $Q$ (up to some tolerance), the K-step becomes $(Q^T A Q) Z + Z M^T = Q^T C V$, where $K = Q Z$. Now the main cost is computing $Q^T A Q$ and $Q^T C V$ because otherwise we have a small Sylvester equation for $Z$. In this handcrafted scenario, the iteration count seen in \Cref{fig:random} is low for both algorithms, since it was made so that the condition number is itself close to 1. Furthermore, both algorithms completely agree on the singular values attained.

\subsubsection{Tucker Tensor Case}

\begin{figure}[ht!]
	\centering
	\includegraphics[width=0.45\linewidth]{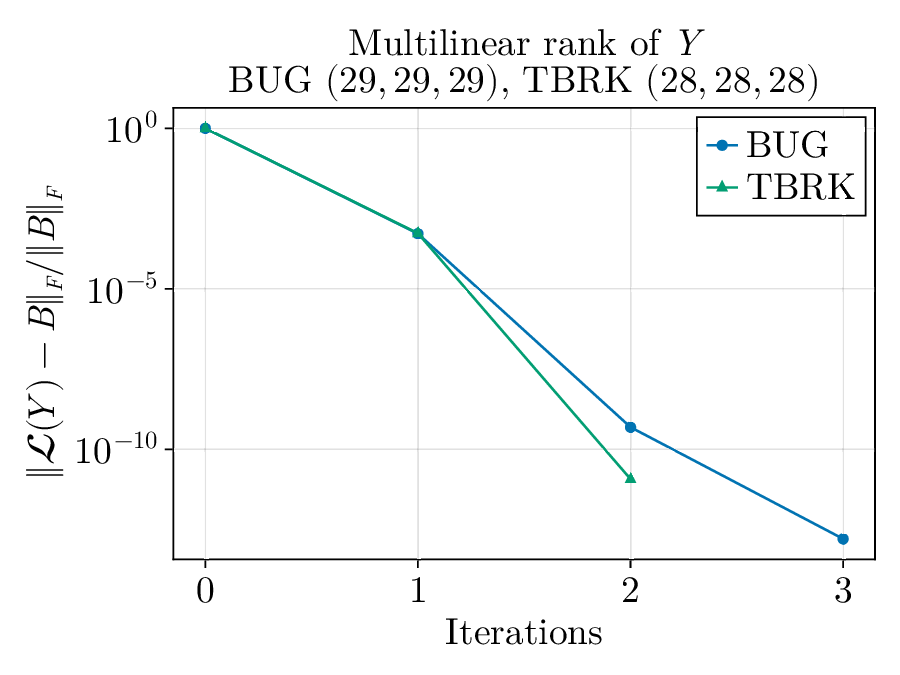}
	\caption{Convergence of the rank-adaptive algorithm and the block rational Krylov method by \cite{Casulli2025} for Tucker tensors with random coefficient matrices of prescribed spectra. The truncation tolerance was set to $\vartheta = 10^{-14}\|C_Y\|_F$.}
	\label{fig:tucker-random}
\end{figure}

Now that we have shown that even $\dist=10$ is enough to provide reasonable low-rank approximations, we use $\dist=100$ in the higher dimensional scenario, for the sake of making computations quicker by conjecturing that the true solution exhibits a stronger ``singular value decay". \Cref{fig:tucker-random} supports this claim as the number of iterations either shrank or stayed the same (though with two iterations the BUG method was close to the threshold of $10^{-10}$ of the relative residual).

\begin{figure}[ht!]
	\centering
	\includegraphics[width=0.9\linewidth]{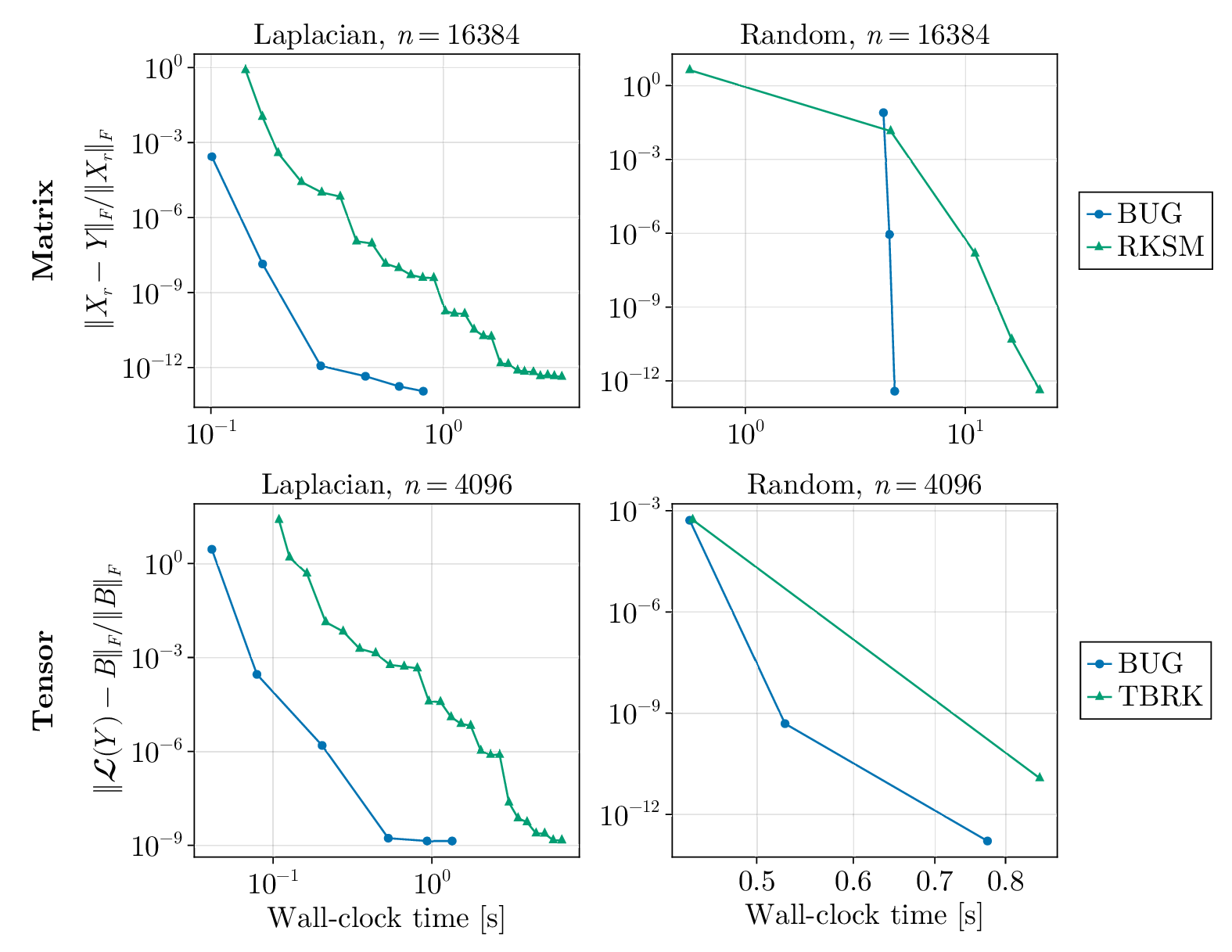}
	\caption{{\bf Top Row}: Forward relative error for the matrix cases against the wall-clock time for both cases shown. {\bf Bottom Row}: Relative Residual for the tensor cases against wall-clock time for both cases.}
	\label{fig:time-error}
\end{figure}

Finally, in \Cref{fig:time-error}, the previous results are shown in terms of wall-clock time using a Intel(R) Xeon(R) CPU E5-2680 v4 @ 2.40GHz.

\section{Discussion}

The experiments with both sparse (Laplacian) and random coefficient matrices demonstrate that the BUG Sylvester solver consistently produces accurate low-rank approximations with few iterations, even in particularly unfavorable scenaria (small eigenvalues of the coefficient matrices). Furthermore, the results for random matrices with prescribed spectra provide empirical evidence that increased spectral separation of the coefficient matrices yields solutions of lower numerical rank, without further conditions on the matrices, such as normality, in accordance to our conjecture. While these findings support the method’s robustness, a rigorous convergence proof remains an important direction for future work.

\section*{Acknowledgements}

I am deeply grateful to Prof. Christian Lubich for encouraging me to work on this method, for the countless in-depth discussions, for helping me with the manuscript and for his integral supervision. I would also like to thank Dr. Dominik Sulz for his help with the enormous amount of questions I had on the generalization of low-rank objects such as the Tucker tensors and further hierarchical structures, and Dr. Yoann Le H\'enaff for valuable discussions in the early stages of this work, when the method was first being developed.

\section*{Declaration of AI Usage}

Artificial Intelligence was only used for the translation of the MATLAB files of the block rational Krylov method of \cite{Casulli2024-vd,Casulli2025} to Julia and the subsequent optimizing, testing and documenting of the proposed method in Julia.

\bibliographystyle{abbrv}
\bibliography{refs}

\end{document}